%% Plain TeX file, the first version was
%% written March-April 1999 during the war against NATO.
%% This is the version of 23 August 1999,  to appear in the
%% "Bulletin" of the Serbian Academy of Sciences.
%% Formula (3.13) corrected December 29, 2000.

\magnification=1200
\nopagenumbers
\font\aa=cmss12

\font\cc=cmcsc8
\font\dd=cmr12
\font\ee=cmr9

\def\txt#1{{\textstyle{#1}}}
\baselineskip=13pt
\def\hf{{\textstyle{1\over2}}}
\def\a{\alpha}\def\b{\beta}
\def\d{{\,\rm d}}
\def\e{\varepsilon}

\def\g{\gamma}
\def\G{\Gamma}
\def\K{{\cal K}}
\def\s{\sigma}

\def\={\;=\;}

\def\zt{\zeta(\hf+it)}

\def\D{\Delta}
  
\def\z{\zeta}
 \def\d{\,{\rm d}} 
\def\R{\Re{\rm e}\,} \def\I{\Im{\rm m}\,} \def\s{\sigma}
\def\={\,=\,}

\def\zt{\zeta(\hf+it)}
\def\Res{\mathop{\rm Res}}
\font\teneufm=eufm10
\font\seveneufm=eufm7
\font\fiveeufm=eufm5
\newfam\eufmfam
\textfont\eufmfam=\teneufm
\scriptfont\eufmfam=\seveneufm
\scriptscriptfont\eufmfam=\fiveeufm
\def\mathfrak#1{{\fam\eufmfam\relax#1}}

\font\tenmsb=msbm10
\font\sevenmsb=msbm7
\font\fivemsb=msbm5
\newfam\msbfam
     \textfont\msbfam=\tenmsb
      \scriptfont\msbfam=\sevenmsb
      \scriptscriptfont\msbfam=\fivemsb
\def\Bbb#1{{\fam\msbfam #1}}

\def \NN {\Bbb N}

  \def\rightheadline{{\hfil{\ee
  On the integral of the error term in the Dirichlet divisor problem}
\hfil\tenrm\folio}}

  \def\leftheadline{{\tenrm\folio\hfil{\ee
  Aleksandar Ivi\'c }\hfil}}
  \def\emptyheadline{\hfil}
  \headline{\ifnum\pageno=1 \emptyheadline\else
  \ifodd\pageno \rightheadline \else \leftheadline\fi\fi}

  \centerline{\aa  ON THE INTEGRAL OF THE ERROR TERM}

\medskip
\centerline{\aa IN THE DIRICHLET DIVISOR PROBLEM}
  \bigskip

\bigskip\bigskip\centerline{\dd  A. IVI\'C }

\bigskip\bigskip
\centerline {(Presented at the 7th Meeting, held on October 29, 1999)}
\bigskip\bigskip
{\ee A b s t r a c t. Several results are obtained concerning
the function $\D_k(x)$, which represents the error term in
the general Dirichlet divisor problem. These include the
estimates for the integral of this function, as well as
for the corresponding mean square integral. The mean square
integral of $\D_2(x)$ is investigated in detail.


}
\bigskip\bigskip

\def\M{{\cal M}}
\def\DJ{\leavevmode\setbox0=\hbox{D}\kern0pt\rlap
 {\kern.04em\raise.188\ht0\hbox{-}}D}

\centerline{\dd 1. Introduction}

\bigskip
Let as usual $d_k(n)\,$ ($k \in \NN$) denote the
number of ways $n$ may be written as a product of $k$ fixed factors. Thus
$d_1(n) \equiv 1, d_2(n) = d(n) = \sum_{\delta|n}1$ is the number of divisors
of $n$, and in general
$$
\sum_{n=1}^\infty\,d_k(n)n^{-s} \;=\;\z^k(s)\qquad(s = \s + it,\,
\s = \R s > 1),\eqno(1.1)
$$
where $\z(s)$ is the Riemann zeta-function. The {\it general (Dirichlet)
divisor problem} (or the Piltz divisor problem, as it is also sometimes
called) consists of the estimation of the quantity
$$
\D_k(x) \;:=\; {\sum_{n\le x}}'d_k(n) - xP_{k-1}(\log x) - \z^k(0),
\eqno(1.2)
$$
which represents the error term in the asymptotic formula for the
summatory function of $d_k(n)$,
where ${\sum\limits_{n\le x}}'$ in general means that the last term in the sum is to
be halved if $x \in\,\NN$, and $\z(0) = -\hf$.
The main term in the formula is
$$
xP_{k-1}(\log x) \;=\; \Res_{s=1}\,\z^k(s)x^ss^{-1},\eqno(1.3)
$$
where $P_{k-1}(u)$ is a certain polynomial in $u$ of degree $k - 1$
whose coefficients depend on
$k$, and $\Res_{z=z_0}F(z)$ denotes the residue of $F(z)$ at the pole
$z = z_0$.
 One has, for example, $P_2(u) = u + 2\gamma - 1$,
where $\gamma = 0.577\ldots\,$ is Euler's constant.
Sometimes the function $\D_k(x)$ is not defined by (1.2) but by
$$
\D_k(x) \;:=\; \sum_{n\le x }d_k(n) - xP_{k-1}(\log x),
$$
but we shall adhere to (1.2).
In general, the coefficients of $P_{k-1}(u)$ may be
found by the use of the Laurent expansion of $\z(s)$ near $s = 1$, namely
$$
\z(s) \;=\; {1\over s - 1} + \sum_{k=0}^\infty \gamma_k(s - 1)^k,
$$
where $(\gamma \equiv \gamma_0)$ for $k \ge 0$
$$
\gamma_k \;=\; {(-1)^k\over k!}\,\lim_{N\to\infty}\left(\sum_{n\le N}
{\log^kn\over n} - {\log^{k+1}N\over k + 1}\right)
$$
are the so-called {\it Stieltjes constants}.  The coefficients of
$P_{k-1}(u)$ in (1.3) were evaluated explicitly in term of the
$\gamma_k$'s by A.F. Lavrik [7].

\smallskip
A large literature exists on estimates for $\D_k(x)$, both pointwise
and in the mean square sense (see e.g.,
[3, Chapters 3, 13] and [12, Chapter 12]).
Of special interest is the function $\D_2(x) \equiv \D(x)$, which represents
the error term in the classical Dirichlet divisor problem, and which
admits an explicit formula, due to F.G. Voronoi (see e.g., [3, Chapter 3]),
in term of the Bessel functions.  The aim of this paper is to provide
estimates for $\int_1^x\D_k(u)\d u$, both pointwise and in the mean
square sense, for certain ``small" values of $k$. This problem appears
to be really of interest for $k > 4$. Namely from Voronoi's formula one has
($c_1$ is a constant)
$$
\eqalign{\int_1^x\D(u)\d u &\=
{1\over2\sqrt{2}\pi^2}x^{3/4}\sum_{n=1}^\infty d(n)n^{-5/4}\sin(4\pi
\sqrt{nx} - {\pi\over4})\cr&
+ {15\over2^6\sqrt{2}\pi^3}x^{1/4}\sum_{n=1}^\infty d(n)n^{-7/4}
\cos(4\pi\sqrt{nx} - {\pi\over4}) + c_1 + O(x^{-1/4}).\cr}\eqno(1.4)
$$
The series in (1.4) are both absolutely convergent and (see e.g., [4])
and they are both also $\Omega_\pm(1)$. Moreover, from (1.4) one
easily obtains
$$
\int_1^X\left(\int_1^x\D(u)\d u\right)^2\d x \;\sim\; CX^{5/2}
\qquad(C > 0,\,X\to\infty).
$$
This settles the case $k = 2$, and for $k = 3$ note that by the complex
integration method (the
Perron inversion formula) one has (see [12, Chapter 12]),
for $X \le x \le 2X$,
$$
\D_3(x) = {1\over\pi\sqrt{3}}x^{1/3}\sum_{n\le X^2}d_3(n)n^{-2/3}
\cos\left(6\pi(nx)^{1/3}\right) + O_\e(X^\e).\eqno(1.5)
$$
Here, as usual, $\e$ denotes arbitrarily small constants which are not
necessarily the same ones at each occurrence.
From (1.5) we obtain
$$
\eqalign{
\int_X^{2X}\D_3(x)\d x &\=
{1\over\pi\sqrt{3}}\sum_{n\le X^2}d_3(n)n^{-2/3}\int_X^{2X}
x^{1/3}\cos\left(6\pi(nx)^{1/3}\right)\d x + O_\e(X^{1+\e})\cr&
\,\ll_\e X\sum_{n\le X^2}d_3(n)n^{-2/3}n^{-1/3} + X^{1+\e} \ll_\e X^{1+\e}
\cr}\eqno(1.6)
$$
by using the first derivative test ([3, Lemma 2.1]). Consequently (1.6)
gives
$$
\int_1^x\D_3(u)\d u \ll_\e x^{1+\e},\quad
\int_1^X\left(\int_1^x\D_3(u)\d u\right)^2\d x \ll_\e X^{3+\e}.\eqno(1.7)
$$
In general, one can obtain by complex integration methods the expression
($1 \ll N \ll x^C$)
$$\eqalign{
\D_k(x) &\=
{x^{{k-1\over2k}}\over\pi\sqrt{k}}\sum_{n\le N}
d_k(n)n^{-{k+1\over2k}}\cos\left(2k\pi(xn)^{{1\over k}} + {(k-3)\pi\over4}
\right) \cr&\,+ \,O_{k,\e}\left(x^\e(1 + x^{{k-1\over k}}N^{-{1\over k}}
+ (xN)^{{1\over2}-{1\over k}})\right).\cr}\eqno(1.8)
$$
However, already for $k = 4$ this formula does not lead to good results.
Namely in evaluating $\int_X^{2X}\D_4(x)\d x$ we shall encounter
$$\eqalign{&
{1\over2\pi}\sum_{n\le N}d_4(n)n^{-5/8}\int_X^{2X} x^{3/8}
\cos\left(8\pi(xn)^{1/4} + {\pi\over4}\right)\d x\cr&
\ll \sum_{n\le N}d_4(n)n^{-5/8}\cdot X^{9/8}n^{-1/4} \ll X^{9/8}
N^{1/8}\log^3N,\cr}
$$
which coupled with the contribution of the error terms in (1.8), will
give a poor final result. For this reason in
the next section we shall adopt another
approach. We shall use power moment results for $\z(s)$  to derive
results on the estimation of $\int_1^x\D_k(u)\d u$. In Section 3
mean square results will be discussed, and in Section 4 we shall discuss
mean square results for $\D_k(x)$, with the accent on the most
important case $k = 2$.

\bigskip
\centerline{\dd 2. Estimates for the integral of the error term}

\bigskip

We start from the classical Perron inversion formula which gives, for
suitable $0 < c < 1$,
$$
\D_k(x) \= {1\over2\pi i}\int_{(c)}\z^k(s)\,{x^s\over s}\d s + (-\hf)^k,
\eqno(2.1)
$$
where as usual
$$
\int_{(c)} F(s)\d s \= \lim_{T\to\infty}\int_{c-iT}^{c+iT} F(s)\d s.
$$
Integration of (2.1) gives then, for $x > 1$,
$$\eqalign{
\int_1^x \D_k(u)\d u &\=
{1\over2\pi i}\int_{(c)}\z^k(s)\,{x^{s+1}\over s(s+1)}\d s + O(x)\cr&
\ll x^{c+1}\int_{-\infty}^{\infty}\,{|\z(c+it)|^k\over1+t^2}\d t
+ x \ll x^{c+1}\cr}\eqno(2.2)
$$
with $c = \eta_k + \e$, where $\eta_k$ is the infimum of $\eta$ for which
one has
$$
\int_T^{2T}|\z(\eta+it)|^k\d t \;\ll_\e\; T^{2+\e}.\eqno(2.3)
$$
We shall obtain
$$
\eta_k \;\le\;{1\over2} - {1\over k}\qquad(4 \le k \le 8).\eqno(2.4)
$$
The proof of (2.4) will be given now in the most interesting case
$k = 8$. By using the functional equation ([3, Chapter 1])
$$
\z(s) \= \chi(s)\z(1-s),\quad\chi(s) \= 2^s\pi^{s-1}\sin({\pi s\over2})
\Gamma(1-s) \;\asymp\; |t|^{{1\over2}-\s}
\;\;(s = \s + it)
$$
we have, for $0 < \eta < 1$,
$$
\int_T^{2T}|\z(\eta+it)|^8\d t \;\ll\; T^{4-8\eta}\int_T^{2T}|\z(1-\eta+it)|^8\d t.
\eqno(2.5)
$$
Now we take $\eta ={3\over8}$ and use the bound ([3, Chapter 8])
$$
\int_T^{2T}|\z({\txt{5\over8}}+it)|^8\d t \;\ll_\e\; T^{1+\e}\eqno(2.6)
$$
to obtain from (2.5)
$$
\int_T^{2T}|\z({\txt{3\over8}}+it)|^8\d t
\ll T\int_T^{2T}|\z({\txt{5\over8}}+it)|^8\d t  \ll_\e T^{2+\e},
$$
which gives $\eta_8 \le {3\over8}$, as asserted. One actually has
$\eta_k = \hf - {1\over k}$ for $2 \le k \le 8$, which is not difficult to
see. Unfortunately the existing results on power moments of $\z(s)$ do
not permit one to extend the validity of (2.4) to any $k$ satisfying
$k > 8$. Nevertheless one can find an upper bound for $\eta_k$ for
any given $k > 8$, but a general expression for this upper
bound, considered as a function of $k$, would be rather complicated.
For this reason we shall content ourselves with explicit bounds for
``small" values of $k$, in particular for $k \le 12$.

From the bounds (see [3, Chapter 8])
$$
\int_0^T|\z(\hf + it)|^9\d t \ll_\e T^{{13\over8}+\e},
\quad \int_0^T|\z({\txt{35\over54}} + it)|^9\d t \ll_\e T^{1+\e}
$$
and convexity of mean values ([3, Lemma 8.3]) one obtains
$$
\int_0^T|\z(\s + it)|^9\d t \ll_\e T^{{239-270\s\over64}+\e}
\qquad(\hf \le \s \le {\txt{35\over54}}).
$$
This yields, for ${19\over54} \le \eta \le \hf$,
$$\eqalign{&
\int_T^{2T}|\z(\eta + it)|^9\d t \ll T^{9({1\over2}-\eta)}
\int_T^{2T}|\z(1-\eta + it)|^9\d t\cr&
\ll_\e T^{9({1\over2}-\eta)+{270(\eta-1)+239\over64}+\e} \ll_\e T^{2+\e}\cr}
$$
for $129 \le 306\eta$, giving
$$
\eta_9 \;\le\; {43\over102} \= 0.421568627\ldots\;.
$$

For the case $k = 10$ we use the bound  of Zhang [13] (the bound
of Ivi\'c--Ouellet [6, p. 250] is slightly weaker, leading to
$\eta_{10} \le {73\over160} = 0.45625$)
$$
\int_T^{2T}|\z(\s + it)|^{10}\d t \ll_\e T^{{17-20\s\over4}+\e}
\qquad({\txt{9\over20}} \le \s \le \hf)
$$
to obtain that $(17-20\s)/4 \le 2$ for $\s \ge 9/20$, giving
$$
\eta_{10} \;\le\; {9\over20} \= 0.45.
$$
Similarly from the bounds ([3, Chapter 8])
$$
\int_0^T|\z(\hf + it)|^{11}\d t \ll_\e T^{{15\over8}+\e},
\quad
\int_0^T|\z({\txt{7\over10}} + it)|^{11}\d t \ll_\e T^{1+\e}
$$
we obtain
$$
\eta_{11} \;\le\; {51\over106} \= 0.481132075.
$$
The slightly better bound
$$
\int_0^T|\z({\txt{1232\over1771}} + it)|^{11}\d t \ll_\e T^{1+\e},
\quad {1232\over1771} = 0.6956521\ldots\;,
$$
of [6] would give a further slight improvement of the bound for
$\eta_{11}$. Finally from
$$
\int_0^T|\z(\hf + it)|^{12}\d t \ll_\e T^{2+\e}
$$
it follows that $\eta_{12} \le \hf$.

\smallskip
Now let $\theta_k$ denote the infimum of $\theta \,(> 0)$ for which
$$
\int_1^x\D_k(u)\d u \;\ll\;x^\theta,\eqno(2.7)
$$
and as usual let $\b_k$ denote the infimum of $b$ for which
$$
\int_1^x\D_k^2(u)\d u \;\ll\; x^{1+2b}.\eqno(2.8)
$$
Then we have
$$
\theta_k \le 1 + \eta_k\,\qquad(k \ge 4),\eqno(2.9)
$$
and from [3, Lemma 13.1] (with $\eta_{2k} = \g_k$) we have $\eta_{2k} = \b_k$,
hence for even $k$ we can bound $\theta_k$ in terms of $\b_k$ and (2.9).
Collecting the above results we obtain

\bigskip
THEOREM 1. {\it We have the bounds
$$
\theta_3 \le 1,\, \theta_k \le {3\over2} - {1\over k} \;(4 \le k \le 8),\,
\theta_9 \le {145\over102}, \,\theta_{10} \le {29\over20}\,,
\,\theta_{11} \le {157\over106},\, \theta_{12}
\le {3\over2},
$$
and in general for $k \ge 2$ we have}
$$
\theta_{2k} \;\le\; 1 + \b_k.\eqno(2.10)
$$
\medskip

 It is known that $\b_k \ge (k-1)/(2k)$
for $k \ge 2$, and in fact the Lindel\"of hypothesis ($\zt \ll_\e |t|^\e$)
is equivalent to $\b_k \= (k-1)/(2k)$ for every $k \ge 2$. We
know at present that  $\b_k \= (k-1)/(2k)$ holds for $k = 2,3,4$, while
e.g. $\b_5 \le {9\over20}$ (see [13]) and
$\b_6 \le \hf $ (see [3]).
It is not easy to surmise what is the true value of $\theta_k$,
and in particular to see how sharp is the inequality in (2.9).

\bigskip
\centerline{\dd 3. The mean square of the integral of the error term}

\bigskip

The approach to mean square estimates for $\int_1^x\D_k(u)\d u$ is
based on the use of Parseval's formula for Mellin transforms (see
E. C. Titchmarsh [11]).
The Mellin transform of an integrable function $f(x)$ is commonly defined
as
$$
\M[f(x)] \, = \, F(s) = \int_0^\infty x^{s-1}f(x)\d x\quad(s = \sigma + it).
$$
An important   feature of Mellin transforms is
the so-called inversion formula. It states that
if $F(s) = \M[f(x)],\;$
$ y^{\sigma-1}f(y) \in L^1(0,\,\infty)$ and $f(y)$
is of bounded variation in a
neighbourhood of $y = x$, then
$$
{f(x + 0) + f(x - 0)\over2} \,
= \, {1\over2\pi i}\int_{(\sigma)}\,F(s)x^{-s}\d s.
\eqno(3.1)
$$
Conversely if (3.1) holds, then $F(s) = \M[f(x)]\,$. We recall that
if $f(x)$ denotes measurable functions, then
$$
L^p(a,b) := \left\{f(x)\Big|\;
\int_a^b|f(x)|^p\d x < +\infty\right\}.
$$
A form of Parseval's formula for Mellin transforms  is the relation
$$
\int_0^\infty f(x)g(x)x^{2\s-1}\d x = {1\over2\pi i}\int_{(\s)}
F(s)\overline{G(s)}\d s,\eqno(3.2)
$$
which holds e.g., if
$$
F(s) = \M[f(x)],\; G(s) = \M[g(x)],\; x^{\s-{1\over2}}f(x)
\in L^2(0,\,\infty),\;x^{\s-{1\over2}}g(x)
\in L^2(0,\,\infty).
$$

\smallskip
The starting point for our mean square results is the bound
$$
\int_X^{2X}\left(\int_1^x\D_k(u)\d u\right)^2\d x \;\ll\;
\int_X^{2X}|f_k(x)|^2\d x + X^3,\eqno(3.3)
$$
which follows from (2.2) and (3.2) with
$$
f_k(x) \;:=\; {1\over2\pi i}\int_{(c)}\,{\z^k(s)x^{s+1}\over s(s+1)}\d s
\qquad(c > \eta_k)\eqno(3.4)
$$
for $x \ge 1$, and $f_k(x) = 0$ for $x < 1$. From (3.1) and (3.4) we have
$$
{\z^k(s)\over s(s+1)} \= \int_0^1 f_k\left({1\over x}\right)x\cdot x^{s-1}\d x
\qquad(\s > \eta_k).\eqno(3.5)
$$
Consequently (3.2) yields
$$
\int_0^1\left|f_k\left({1\over x}
\right)\right|^2x^2x^{2c-1}\d x \= {1\over2\pi}\int_{-\infty}
^\infty \,{|\z(c+it)|^{2k}\over(c^2+t^2)((c+1)^2+t^2)}\d t.\eqno(3.6)
$$
Since, by $L^2$--theory, the convergence of one integral in (3.6) implies
the convergence of the other one, and
$$
\int_0^1\left|f_k\left({1\over x}\right)\right|^2x^2x^{2c-1}\d x \=
\int_1^\infty |f_k(x)|^2 x^{-3-2c}\d x, \eqno(3.7)
$$
it follows from (3.6) and (3.7) that
$$
\int_X^{2X} |f_k(x)|^2\d x \;\ll\; X^{3+2c_k}\qquad(X > 1),\eqno(3.8)
$$
provided that $0 < c_k < 1$ is such a constant for which
$$
\int_{-\infty}^\infty \,|\z(c_k+it)|^{2k}\,{\d t\over1+t^4}
\;\ll\;1.
$$
The last condition reduces to finding $0 < \s_k < 1$ such that
$$
\int_T^{2T} |\z(\s_k + it)|^{2k}\d t \;\ll_\e\; T^{4+\e},\eqno(3.9)
$$
and then one can take $c_k = \s_k + \e$ in (3.8). We trivially have
$c_3 = \e$ (see (1.7)), and also $c_4 = \e$ (follows from $\z(it)
\ll t^{1/2}\log t$).

For $k > 4$ we have, by the functional equation for $\z(s)$,
$$
\int_T^{2T}|\z(\s_k + it)|^{2k}\d t \ll T^{k-2k\s_k}
\int_T^{2T}|\z(1-\s_k + it)|^{2k}\d t \ll_\e T^{1+k-2\s_k+\e},\eqno(3.10)
$$
provided that for a given $k$ one can find $0 < \s_k < 1$ for which one has
$$
\int_0^T|\z(1-\s_k+it)|^{2k}\d t \;\ll_\e\; T^{1+\e}.\eqno(3.11)
$$
If we can take
$$
\s_k \= {k-3\over2k},\eqno(3.12)
$$
then from (3.11) we obtain
$$
\int_T^{2T}|\z(\s_k + it)|^{2k}\d t \;\ll_\e\; T^{4+\e},
$$
which shows that $c_k = (k-3)/(2k) + \e$ is permissible. We can infer
that (3.11) holds with (3.12) if $k=5$ and $k=6$, since we have the
bound ([3, Chapter 8])
$$
\int_0^T|\z({\txt{3\over4}}+it)|^{12}\d t \;\ll_\e\; T^{1+\e}.
$$
It is very likely that $c_k = {k-3\over2k} + \e$ will hold for
at least some $k > 6$, but this cannot be inferred from the existing
results on power moments of $\z(s)$. For $k > 6$ one will have to use
a weaker bound than (3.11), and consequently we shall have a weaker
bound than $c_k = {k-3\over2k} + \e$. A general formula for $c_k$ is
possible, but its form would be rather complicated. For this reason
we shall content ourselves with the above bounds, which we formulate as

\bigskip
THEOREM 2. {\it Let $\rho_k$ be the infimum of $\rho > 0$ for which
$$
\int_1^X\left(\int_1^x\D_k(u)\d u\right)^2\d x \;\ll\; X^\rho.
$$
Then we have} $\rho_3 \le 3, \rho_4 \le {13\over4}, \rho_5 \le {17\over5},
\rho_6 \le {7\over2}$.

\bigskip
From the definition of $\rho_k$ it easily follows that
$$
\rho_k \;\le\; 3 + 2\b_k,\eqno(3.13)
$$
where $\b_k$ is as in Section 2. Note, however, that the bounds of Theorem 3
are much better than the bounds that one can derive from (3.13) and
the sharpest known bounds for $\b_k$.

\bigskip
\centerline{\dd 4. The mean square formula for $\D_k(x)$}

\bigskip
Let us define, for $k \ge 2$,
$$
\K_k(s) \= \int_1^\infty \D_k^2(x)x^{-s}\d x.\eqno(4.1)
$$
By an integration by parts and the use of (3.13)
it follows that $\K_k(s)$ is a regular function
of $s$ for $\s = \R s > 1 + 2\b_k$. The analytic behaviour of $\K_k(s)$
enables one to obtain information on the mean square of $\D_k(x)$
via the formula
$$
\int_1^X \D_k^2(x)\d x \= {1\over2\pi i}\int_{(1+2\b_k+\e)}\K_k(s)
{X^s\over s}\d s\qquad(X > 1).\eqno(4.2)
$$
Namely by using the classical integral $(c > 0)$
$$
{1\over2\pi i}\int_{(c)}{y^s\over s}\d s = \cases{1 \qquad (y > 1),&\cr  &\cr
\hf\qquad(y = 1),\cr &\cr0\qquad(0 < y < 1),\cr}
$$
we have
$$\eqalign{{1\over2\pi i}\int_{(1+2\b_k+\e)}\K_k(s){X^s\over s}\d s&
= \int_1^\infty\left({1\over2\pi i}\int_{(1+2\b_k+\e)}\left({X\over x}
\right)^s{\d s\over s}\right)\D_k^2(x)\d x \cr&
= \int_1^X \D_k^2(x)\d x\qquad(X > 1).\cr}
$$
We can obtain analytic continuation of $\K_k(s)$ to the left of the
line $\s = 1 + 2\b_k$ in two cases: $k = 2$ and $k = 3$, which follows
from mean square results on $\D_k(x)$ (see [3, Chapter 13]). In the latter
case we use the asymptotic formula
$$
\eqalign{&
\int_1^x\D_3^2(y)\d y \= Cx^{5/3} + R(x),\cr&
C = {1\over10\pi^2}\sum_{n=1}^\infty d_3^2(n)n^{-4/3},\quad R(x)
\ll_\e x^{{14\over9}+\e}.\cr}\eqno(4.3)
$$
From (4.3) we obtain
$$\eqalign{
\K_3(s) &\= \int_1^\infty \D_3^2(x)x^{-s}\d x \= \int_1^\infty
\left({\txt{5\over3}}Cx^{2/3} + R'(x)\right)x^{-s}\d x\cr&
\= {5C\over3s-5} + C_1 + s\int_1^\infty R(x)x^{-s-1}\d x.\cr}\eqno(4.4)
$$
The formula (4.4) holds initially for $\s > 5/3$, but the upper bound
for $R(x)$ in (4.3) shows that it provides analytic continuation
of $\K_3(s)$ to the half-plane $\s > 14/9$, where $\K_3(s)$ is regular
except for a simple pole at $s = 5/3$.

\medskip
The case $k =2$ is even more interesting. We have
$$
\eqalign{&
\int_1^x\D^2(y)\d y \= Ax^{3/2} + F(x),\cr&
A = {1\over6\pi^2}\sum_{n=1}^\infty d^2(n)n^{-3/2},\quad F(x)
\ll x\log^4x,\cr}\eqno(4.5)
$$
where the upper bound for $F(x)$ is due to E. Preissmann [9].
Similarly to (4.4) we obtain from (4.5)
$$
\K_2(s) \= {3A\over2s-3} + C_2 + s\int_1^\infty F(x)x^{-s-1}\d x.\eqno(4.6)
$$
In view of the upper bound for $F(x)$ in (4.5) it follows that (4.6)
provides analytic continuation of $\K_2(s)$ to the half-plane $\s > 1$,
where $\K_2(s)$ is regular, except for a simple pole at $s = 3/2$.
In general, it seems very likely that $\K_k(s)$ possesses analytic
continuation to the left of $\s = 1 + 2\b_k$, and that it has a simple
pole at $s = 1 + 2\b_k$. However, the existing mean square results on
$\D_k(x)$ are not sharp enough to deduce this assertion. In the case of
$\K_2(s)$ one can obtain analytic continuation of $\K_2(s)$ to the
half-plane $\s > 2/3$, as well as mean square results for $\s > 1$.
The results are contained in

\bigskip
THEOREM 3. {\it The function $\K_2(s)$ possesses analytic continuation
to the region $\s > 2/3$, where it is a regular function of $s$, except
at $s = 3/2$ where it has a simple pole, and at $s = 1$ where it has
a pole of order} 3. {\it
In the region $\s > 1$, without an $\e$--neighbourhood
of $s = 3/2$, it is of polynomial growth in $|\I s|$. Moreover,}
$$
\int_0^T|\K_2(\s + it)|^2\d t \;\ll\; \cases{T^{6-4\s}(\log T)^{12-8\s}
&\qquad$(1 < \s < 3/2)$,\cr\cr
1 &\qquad$(\s > 3/2)$.\cr}\eqno(4.7)
$$

\bigskip
{\bf Proof.} We use the Laplace transform formula
$$
\int_0^\infty \D^2(x)e^{-x/T}\d x = {B\over8}\left({T\over\pi}\right)^{3/2}
+ (A_1\log^2T + A_2\log T + A_3)T + O_\e(T^{2/3+\e}),\eqno(4.8)
$$
where
$$
B \= \sum_{n=1}^\infty d^2(n)n^{-3/2},\quad A_1 = -{1\over4\pi^2},
$$
which was proved in [5]. Since the integral defining $\K_2(s)$ is
absolutely convergent for $\s > 3/2$, it follows that for $c > 3/2$ and
$T > 0$ one has
$$\eqalign{&
{1\over2\pi i}\int_{(c)}\G(s)T^s\K_2(s)\d s = \int_1^\infty\D^2(x)
\left({1\over2\pi i}\int_{(c)}\left({x\over T}\right)^s\G(s)\d s\right)\d x
\cr& = \int_1^\infty \D^2(x)e^{-x/T}\d x = \int_0^\infty \D^2(x)e^{-x/T}\d x
+ O(1).\cr}\eqno(4.9)
$$
Hence (4.8) gives, for $c > 3/2$ and $T \ge 1$,
$$
{1\over2\pi i}\int_{(c)}\G(s)T^s\K_2(s)\d s =
{B\over8}\left({T\over\pi}\right)^{3/2}
+ (A_1\log^2T + A_2\log T + A_3)T + O_\e(T^{2/3+\e}),
$$
or for $0 < x \le 1, c > 3/2$,
$$
\eqalign{f(x) \;&:=\;
{1\over2\pi i}\int_{(c)}\G(s)x^{-s}\K_2(s)\d s \cr&\;=
{B\over8}(\pi x)^{-3/2}
+ (A_1\log^2({1\over x}) + A_2\log ({1\over x}) + A_3){1\over x}
+ O_\e(x^{-2/3-\e}).\cr}\eqno(4.10)
$$
From (4.10) we deduce for $\s > 3/2$, by the Mellin inversion formula
(3.1),
$$
\G(s)\K_2(s) = \int_0^\infty f(x)x^{s-1}\d x =
\int_0^1 f(x)x^{s-1}\d x + \int_1^\infty f(x)x^{s-1}\d x = I_1(s) + I_2(s),
\eqno(4.11)
$$
say. Note that the definition of $f(x)$ in (4.10) makes sense
if $x > 0$, and (4.9) with $T = 1/x$ yields
$$
f(x) \;\ll\;e^{-{1\over2}x}\qquad(x \ge 1),\eqno(4.12)
$$
hence (4.12) shows that $I_2(s)$ is regular for $\s > 0$. To investigate
$I_1(s)$ we use (4.10) to deduce that
$$
I_1(s) = {c_1\over s-{3\over2}} + {a_3\over(s-1)^3} +
{a_2\over(s-1)^2} + {a_1\over s-1} + \int_0^1 h(x)x^{s-1}\d x,\eqno(4.13)
$$
where $c_1,\,a_j$ are effectively computable constants, and
$$
h(x) \;\ll_\e\; x^{-{2\over3}-\e}\qquad(0 < x \le 1).
$$
This means that the integral in (4.13) is regular for $\s > 2/3$,
proving the first part of Theorem 3. From (4.6) we obtain
$$
\K_2(s) \;\ll_\e\; |t|\qquad (\s > 1,\,|s - {\txt{3\over2}}| \ge \e).
$$
It is very likely that the error term in (4.8) can be sharpened to
$O_\e(T^{{1/2}+\e})$. For this reason it also
seems likely that $\K_2(s)$ admits analytic continuation even to
the half-plane $\s > \hf$, where it is of polynomial growth. However
at present it does not seem possible to prove this assertion.

It remains to prove the mean square bounds of (4.7). From Parseval's
formula (3.2) one has, if $\s$ is sufficiently large,
$$
\int_1^\infty \D^4(x)x^{1-2\s}\d x \=
{1\over2\pi}\int_{-\infty}^\infty |\K_2(\s + it)|^2\d t.\eqno(4.14)
$$
Namely (4.14) follows from (4.1) and
$$
\int_1^\infty f^2(x)x^{1-2\s}\d x = {1\over2\pi}
\int_{-\infty}^\infty|F^*(\s + it)|^2\d t,\eqno(4.15)
$$
where
$$
F^*(s) \;:=\; \int_1^\infty f(x)x^{-s}\d x.\eqno(4.16)
$$
One obtains (4.15) from Parseval's formula (3.2) on replacing $f(x)$ and
$g(x)$ by ${1\over x}\bar{f}(x)$, where $\bar{f}(x) = f({1\over x})$
if $0 < x \le 1$ and $\bar{f}(x) = 0$ otherwise.

Note that one has (see D.R. Heath-Brown [2])
$$
\int_0^x\D^4(y)\d y \;\sim\; Cx^2\qquad(C > 0),\eqno(4.17)
$$
hence the integral on the left-hand side of (4.14) is convergent
for $\s > 3/2$. This implies that the integral on the right-hand
side of (4.14) is also convergent
for $\s > 3/2$, giving the second mean square bound in (4.7).

To obtain the first mean square bound in (4.7) write
$$
\K_2(s) = \int_1^\infty \D^2(x)x^{-s}\d x = \int_1^X \D^2(x)x^{-s}\d x +
\int_X^\infty \D^2(x)x^{-s}\d x,
$$
where $X$ will be suitably chosen a little later. Using (4.5) we have
$$\eqalign{
\int_X^\infty \D^2(x)x^{-s}\d x &= \int_X^\infty\left({\txt{3\over2}}Ax^{1/2}
+ F'(x)\right)x^{-s}\d x = {3AX^{{3\over2}-s}\over 2s-3}\cr&
+ O(X^{1-\s}\log^4X) + s\int_X^\infty F(x)x^{-s-1}\d x,\cr}\eqno(4.18)
$$
which provides then the analytic continuation of $\K_2(s)$ to $\s > 1$.
To treat the mean square integral of $\K_2(s)$ when $1 < \s < {3\over2}$
we use the following method. Let us consider
$$
I := \int_T^{2T}\left|\int_a^b g(x)x^{-s}\d x\right|^2\d t
\quad(s = \s + it,\,T \ge T_0 > 0,\, a \ge 1),
$$
and set in (4.15) $f(x) = g(x)$ if $a \le x \le b$ and $f(x) = 0$
otherwise. Then $F^*(s)$ in (4.16) becomes
$$
F^*(s) \= \int_a^b g(x)x^{-s}\d x.
$$
Consequently (4.15) (with $f \equiv g$) gives
$$
{1\over2\pi}I \;\le\; {1\over2\pi}\int_{-\infty}^\infty |F^*(\s + it)|^2\d t
= \int_a^b g^2(x)x^{1-2\s}\d x.\eqno(4.19)
$$
Returning to the mean square of $\K_2(s)$ we have from (4.18)
and (4.19), when
$1 < \s < {3\over2}$,
$$\eqalign{&
\int_T^{2T}|\K_2(\s + it)|^2\d t \ll \int_T^{2T}\left|\int_1^X\D^2(x)x^{-s}
\d x\right|^2\d t +
\int_T^{2T}\left|\int_X^\infty\D^2(x)x^{-s}\d x\right|^2\d t \cr&
\ll \int_1^X \D^4(x)x^{1-2\s}\d x + T^{-1}X^{3-2\s} + TX^{2-2\s}\log^8X
+ T^2\int_X^\infty F^2(x)x^{-1-2\s}\d x\cr&
\ll X^{3-2\s} + TX^{2-2\s}\log^8X + T^2X^{2-2\s}\log^4X \ll
T^{6-4\s}(\log T)^{12-8\s}\cr}
$$
with the choice $X = T^2\log^4T$. Here we used  (4.17) and the
bound of K.-M. Tsang  [10] (with $r=2$)
$$
\int_2^X|F(x) + (4\pi^2)^{-1}x\log^2x - \kappa x\log x|^r\d x \ll (cr)^{4r}
X^{r+1},\eqno(4.20)
$$
which is valid uniformly for $X > 2,\,r \in \NN$ and
suitable constants $\kappa$
and $c$. This completes the proof of Theorem 3.

\medskip
In concluding, it seems in order to discuss the shape of the mean square
formula (4.5). In [8] and [10] some remarkable results on $F(x)$ were
proved by Lau and Tsang, which include the bound (4.20). From (4.20)
Tsang deduces that, for almost all $x$,
$$
F(x) \= -{1\over4\pi^2}x\log^2x + \kappa x\log x + O(x),
\eqno(4.21)
$$
and conjectures that (4.21) holds for all $x \,(\ge 2)$. In view of Theorem
3 and (4.2) it seems that perhaps a formula sharper than (4.21) holds
for all $x$, namely
$$
F(x) \= -{1\over4\pi^2}x\log^2x + \kappa x\log x + \lambda x + G(x),
\;G(x) = O(x^\a),\eqno(4.22)
$$
for suitable constants $\kappa,\,\lambda$ and $0 < \a < 1$. Namely
heuristically we use (4.2) with $k=2$ and shift the line of integration
to the left, passing over the poles at $s = 3/2$ (which yields the main
term $Ax^{3/2}$ in (4.5) and at $s=1$ (which yields the main term
in (4.22)). I conjecture that (4.22) holds with any $\a$ satisfying
${3\over4} < \a < 1$. The reason for the bound $\a > {3\over4}$ is that
(4.22) with $\a < {3\over4}$ is not possible, which will be shown now.
We start from ([6, Lemma 2])
$$
\D(x) \= H^{-1}\int_x^{x+H}\D(y)\d y + O(H\log x)\qquad(x^\e \le H \le x).
\eqno(4.23)
$$
By the Cauchy-Schwarz inequality for integrals and (4.5), (4.23) implies
$$\eqalign{
\D^2(x) &\ll H^{-1}\int_x^{x+H}\D^2(y)\d y  + H^2\log^2x\cr&
\ll x^{1/2} + H^{-1}(F(x+H) - F(x-H)) + H^2\log^2x.\cr}\eqno(4.24)
$$
Now suppose that (4.22) holds with some $\a < {3\over4}$. Since we have
(see e.g., J.L. Hafner [1])
$$
\limsup_{x\to\infty}\,|\D(x)x^{-1/4}| \= \infty,
$$
this means that, given any constant $C > 0$, there exist arbitrarily
large values $X_1$ such that $\D^2(X_1) > CX_1^{1/2}$. Therefore (4.24)
yields
$$\eqalign{
X_1^{1/2} &\ll H^{-1}\Bigl(F(X_1+H) - F(X_1-H)\Bigr) + H^2\log^2X_1\cr&
\ll H^{-1}X_1^\a + H^2\log^2X_1 \ll X_1^{2\a/3}\log^2X_1\cr}\eqno(4.25)
$$
if $H = X_1^{\a/3}$. But since $\a < {3\over4}$, (4.25) gives then
a contradiction, proving the assertion. On the other hand, the above
proof shows that if (4.22) holds with $\a = {3\over4} + \e$, then (4.24)
yields the conjectural bound $\D(x) \ll_\e x^{1/4+\e}$, which is very
strong.

\bigskip
\centerline{\dd REFERENCES}

\bigskip
\item{[1]} J.L. Hafner, New omega theorems for two classical lattice
point problems, Invent. Math. {\bf63}(1981), 181-186.

\item{[2]} D.R. Heath-Brown, The distribution and moments of the error
term in the Dirichlet divisor problem, Acta Arith. {\bf60}(1992), 389-415.

\item{[3]} A. Ivi\'c, The Riemann zeta-function, John Wiley \& Sons,
New York, 1985.

\item{[4]} A. Ivi\'c, Large values of certain number-theoretic error
terms, Acta Arith. {\bf56}(1990), 135-159.

\item{[5]} A. Ivi\'c, The Laplace transform of the square in the circle
and divisor problems, Studia Scient. Math. Hung. {\bf32}(1996), 181-205.

\item{[6]} A. Ivi\'c and M. Ouellet, Some new estimates in the Dirichlet
divisor problem, Acta Arith. {\bf52}(1989), 241-253.

\item{[7]} A.F. Lavrik, On the principal term in the divisor problem and
the power series of the Riemann zeta-function in a neighbourhood of its
pole (Russian), Trudy Mat. Inst. Steklova  {\bf142}(1976), 165-173.

\item{[8]} Y.-K. Lau and K.-M. Tsang, Mean square of the remainder term
in the Dirichlet divisor problem, J. Th\'eorie des Nombres Bordeaux
{\bf7}(1995), 75-92.

\item{[9]} E. Preissmann, Sur la moyenne quadratique du terme de reste
du probl\`eme du cercle, C.R. Acad. Sciences Paris S\'erie I {\bf306}(1988),
151-154.

\item{[10]} K.-M. Tsang, Mean square of the remainder term in the Dirichlet
divisor problem II, Acta Arith. {\bf71}(1995), 279-299.

\item{[11]} E.C. Titchmarsh, Introduction to the theory of Fourier integrals,
Clarendon Press, Oxford, 1948.

\item{[12]} E.C. Titchmarsh, The theory of the Riemann zeta-function
(2nd ed.), Clarendon Press, Oxford, 1986.

\item{[13]} Wengpeng Zhang, On the divisor problem, Kexue Tongbao
{\bf33}(1988), 1484-1485.

\bigskip\bigskip
\parindent=0pt
\cc
Aleksandar Ivi\'c \par
Katedra Matematike RGF-a\par
Universiteta u Beogradu\par
Dju\v sina 7, 11000 Beograd,
\par
Serbia (Yugoslavia)\par
{\sevenbf e-mail: ivic@rgf.bg.ac.yu}

\bye